\documentclass[a4paper]{amsart}
\usepackage[utf8]{inputenc}
\usepackage[T1]{fontenc}
\usepackage[colorinlistoftodos]{todonotes}
\usepackage{xifthen, xfrac, mathrsfs, url, amssymb, enumitem}
\usepackage{hyperref}
\newtheorem{thm}{theorem}
\newtheorem{alg}[thm]{Algorithm}
\theoremstyle{remark}
\newtheorem{rem}{Remark}
\newenvironment{poc}[1][]{\begin{proof}[\ifthenelse{\equal{#1}{}}{Proof of correctness}{Proof of correctness of #1}]}{\end{proof}}
\newcommand{\FF}{\mathbb{F}}
\newcommand{\QQ}{\mathbb{Q}}
\newcommand{\ZZ}{\mathbb{Z}}
\newcommand{\SB}{\mathscr{B}}
\newcommand{\gd}{\mathfrak{d}}
\newcommand{\gdi}{{\gd_i}}
\newcommand{\gp}{\mathfrak{p}}
\newcommand{\gpi}{{\gp_i}}
\newcommand{\gq}{\mathfrak{q}}
\newcommand{\gr}{\mathfrak{r}}
\newcommand{\q}{q}
\newcommand{\K}{K}
\renewcommand{\L}{L}
\newcommand{\Kdi}{K_\gdi}
\newcommand{\Kp}{K_\gp}
\newcommand{\Kpi}{K_\gpi}
\newcommand{\Kq}{K_\gq}
\newcommand{\Kr}{K_\gr}
\newcommand{\card}[1]{|#1|}
\newcommand{\len}{\ell}
\newcommand{\ext}[2]{#1/#2}
\newcommand{\form}[1]{\langle#1\rangle}
\newcommand{\units}[1]{#1^\times}
\newcommand{\un}[1][]{\ifthenelse{\equal{#1}{}}{\units{\K}}{\units{#1}}}
\newcommand{\squares}[1]{#1^{\times2}}
\newcommand{\sq}[1][]{\ifthenelse{\equal{#1}{}}{\squares{\K}}{\squares{#1}}}
\newcommand{\norm}[2][\ext{\L}{\K}]{\ifthenelse{\equal{#2}{}}{N_{#1}}{N_{#1}\left(#2\right)}}

\newcommand{\odd}{\equiv 1\pmod{2}}
\newcommand{\Sing}[1]{\mathbb{E}_{#1}}
\newcommand{\SingS}{\Sing{S}}
\DeclareMathOperator{\ord}{ord}
\DeclareMathOperator{\sgn}{sgn}
\newcommand{\term}[1]{\emph{#1}}
\makeatletter
\@namedef{subjclassname@2020}{\textup{2020} Mathematics Subject Classification}
\makeatother
\title{Solving sums of squares in global fields}
\author{Przemys{\l}aw Koprowski}
\email{przemyslaw.koprowski@us.edu.pl}
\address{Institite of Mathematics, University of Silesia, ul. Bankowa 14, Katowice, Poland}
\begin{document}
\maketitle
\begin{abstract}
The problem of writing a totally positive element as a sum of squares has a long history in mathematics, going back to Bachet and Lagrange. While for some specific rings (like integers or polynomials over the rationals), there are known methods for decomposing an element into a sum of squares, in general, for many other important rings and fields, the problem is still widely open. In this paper, we present an explicit algorithm for decomposing an element of an arbitrary global field (either a number field or a global function field) into a sum of squares of minimal length.
\end{abstract}

\section{Introduction}
One of the fundamental results in classical number theory is the famous Lagrange Four Squares Theorem, stating that every positive integer is a sum of four squares. Over centuries, Lagrange theorem has been generalized in numerous directions. A remarkable result by Siegel says that a totally positive element of any global field is a sum of four squares. Hence, from the existential standpoint, the theory over an arbitrary number field is as complete as over the rationals. It is no longer so if we consider the computational aspect of the problem. There is a number of algorithms for decomposing a positive rational number into a sum of four (or less) squares, see for example \cite{Bumby1996, PS2019, PT2018, RS1986, Schoof1985}. On the other hand, for number fields, or more generally global fields, there is only the recent paper by M.-K.~Darkey-Mensah and B.~Rothkegel \cite{DMR202x}, where the authors give an algorithm for computing the \term{length} (i.e. the minimal number of summands needed to represent given $a$ as a sum of squares) of an ele­ment of a global field. The author is not aware of any algorithm that would actually solve the problem of representing an algebraic number as a sum of squares of minimal length. Thus, the goal of this paper is to fill this gap and propose such algorithms, working over arbitrary global fields (so either number fields or global function fields).

Recall that over a field of characteristic~$2$, every sum of squares is a square. Hence the problem of decomposing an element into a sum of squares is meaningful only over fields of characteristic distinct from~$2$. Therefore, in what follows the term \term{global field} will always mean ``global field of characteristic $\neq 2$''. It can be either a number field (then the characteristic is just zero) or an algebraic function field over a finite field~$\FF_{p^n}$, for some odd prime~$p$.

Throughout this paper we use the following notation and terminology. By~$\K$ we always denote a global field and $\un$ stands for the group of nonzero elements. We write $\form{a_1, \dotsc, a_n}$ for the quadratic form $a_1x_1^2 + \dotsb + a_nx_n^2$. The \term{length} of an element $a\in \un$, denoted $\len(a)$, is the minimal number of summands needed to express~$a$ as a sum of squares. If $a$ is not a sum of squares, we put $\len(a) := \infty$. Since a representation of $-1$ as a sum of squares plays a special role, the length of $-1$ is called the \term{level} of~$\K$ and denoted $s(\K)$. In particular, $s(\K) = \infty$ if and only if $\K$ is formally real. It is known (see e.g., \cite[Example~XI.2.4\;(8)]{Lam2005}) that the level of a global field can take only one of the following values: $1$, $2$, $4$ or~$\infty$. Another invariant related to sums of squares is the \term{Pythagoras number} of~$\K$, denoted $P(K)$ hereafter. By definition, the Pythagoras number is the supremum of the lengths of all the elements of K that are sums of squares.

An equivalency class of valuations is called a \term{place} of~$\K$. We always use fraktur letters $\gd$, $\gp$, $\gq$, $\gr$ to denote places. To avoid a repeated use of the word ``place'', non-archimedean places will be called also \term{primes} (or \term{finite primes}, when we wish to emphasize the fact that they are not archimedean). The com­pletion of~$\K$ with respect to a place~$\gp$ is denoted~$\Kp$. If $\gp$ is non-archime­dean, then $\ord_\gp: \K \to \ZZ\cup\{\infty\}$ stands for the corresponding discrete valuation on~$\K$. The valuation $\ord_\gp$ induces a map $\sfrac{\un}{\sq}\to \sfrac{\ZZ}{2\ZZ}$ on the square class group of~$\K$. Abusing the notation harmlessly, we denote the latter map by $\ord_\gp$, again.

Let $S$ be a finite set of places of~$\K$, containing all archimedean places. A square-class $\alpha\in \sfrac{\un}{\sq}$ is said to be \term{$S$-singular} if $\ord_\gp\alpha$ is even for every prime $\gp\notin S$. The set of $S$-singular square classes form a subgroup of $\sfrac{\un}{\sq}$, that we denote~$\SingS$. Observe that $\SingS$ is an elementary $2$-group, hence it can be treated as a vector space over~$\FF_2$. We will use this fact extensively.

For a real place~$\gr$, by $\sgn_\gr$ we denote the sign function on~$\K$, associated to the ordering of~$\K$, induced by the unique ordering of the completion~$\Kr$. Given a finite prime~$\gp$ and two elements~$a$ and~$b$ of~$\K$, we write $(a, b)_\gp$ for their Hilbert symbol at~$\gp$. Further, if $\q = \form{a_1, \dotsc, a_n}$ is a quadratic form (with coefficients in~$\K$), then $\q\otimes \Kp$ stands for the associated quadratic form over~$\Kp$. We write $s_\gp\q$ for the Hasse invariant of~$\q$ at~$\gp$. By definition (see e.g., \cite[Definition V.3.17]{Lam2005}):
\[
s_\gp\q := \prod_{i < j} (a_i, a_j)_\gp.
\]

Our algorithms proposed in this paper rely on a number of already existing tools. The three main ingredients are:
{
\setlist[enumerate]{wide=\parindent}
\begin{enumerate}
\item Factorization of (fractional) ideals in global fields. There are a couple of methods for achieving this (see e.g., \cite{Cohen1993, DMK2019}). At the moment of writing the fastest ideal factorization method known to the author is the one described in \cite{GMN2011}.
\item Solution of a norm equation of a form $\norm{x} = b$, where $b\in \un$ and~$\L$ is an extension of~$\K$. Algorithms for solving the norm equations are described in \cite{Cohen2000, FJP1997, FP1983, Garbanati1980, Simon2002}. In fact, in this paper, $\ext{\L}{\K}$ will be always a quadratic extension.
\item Construction of a basis (over~$\FF_2$) of the group~$\SingS$ of $S$-singular square classes, for a given set~$S$. The author is aware of three algorithms for this task. One approach, outlined in Magma manual \cite{CBFS2015}, is to find a set~$S'$ containing~$S$ and such that the $S'$-class number is odd. Then $\Sing{S'}$ coincides with the group of $S'$-units modulo squares. The sought group~$\SingS$ can be subsequently constructed as a subspace of~$\Sing{S'}$. Another method due to A.~Czogała has been recently described in \cite{KR202x}. Yet another algorithm is published in \cite{Koprowski2021}. The last one is probably best suited for our purpose, since it is incremental in nature, and our solution is based on a procedure of gradual enlarging a certain set~$S$.
\end{enumerate}
}

\section{Algorithms}
We may assume that the length of~$a$ is known a priori, since it can be calculated using \cite{DMR202x}. Of course, if the length is $\len(a) = 1$, then $a$ is a square and it suffices to compute its square root in~$\K$. Below we present separate algorithms for all other possible lengths of~$a$. We begin with sums of just two squares.

\begin{alg}
Given an element $a\in\un$ of length $\len(a) = 2$, this algorithm outputs $c_1, c_2\in \un$ such that $a = c_1^2 + c_2^2$.
\begin{enumerate}
\item\label{st:len2:-1} If $-1$ is a square in~$\K$, then output 
\[
c_1 := \frac{a+1}{2},\qquad c_2 := \frac{a-1}{2}\cdot \sqrt{-1}
\]
and quit.
\item Otherwise set $\L := \K\bigl(\sqrt{-1}\bigr)$ and solve the norm equation $\norm{x} = a$. Denote a solution by $c_1 + c_2\cdot \sqrt{-1}$.
\item Output $c_1, c_2$.
\end{enumerate}
\end{alg}

\begin{poc}
If $-1$ is a square in~$\K$, then the correctness of the output of step~\eqref{st:len2:-1} follows by a direct calculation. Otherwise, when $-1$ is not a square, then $\L := \K\bigl(\sqrt{-1}\bigr)$ is a proper extension of~$\K$ and we have
\[
a = c_1^2 + c_2^2 = \norm{c_1 + c_2\cdot \sqrt{-1}}
\]
as desired.
\end{poc}

We may now focus on sums of three squares. First of all, observe that if $\len(a) = 3$ for some $a\in \K$, then the Pythagoras number of~$\K$ must be at least three. It follows from \cite[Theorem XI.5.6]{Lam2005} that the level of~$\K$ is strictly grater than~$1$. Hence it is either~$2$, or~$4$, or infinity if $\K$ is a formally real number field. In particular, $-1\notin \sq$ and the field $\L := \K\bigl(\sqrt{-1}\bigr)$ is always a proper extension of~$\K$. We need to distinguish between two cases. For the sake of clarity of the exposition, we present them as two separate algorithms. The first one is rather trivial.

\begin{alg}
Let~$\K$ be a global field of level $s(K) = 2$. Given an element $a\in \un$ of length $\len(a) = 3$, this algorithm outputs $c_1, c_2, c_3\in \un$ such that $a = c_1^2 + c_2^2 + c_3^2$.
\begin{enumerate}
\item Let $\L := \K\bigl(\sqrt{-1}\bigr)$. 
\item\label{st:len3s:norm} Solve the norm equation $\norm{x} = -1$ and let $d_1 + d_2\cdot \sqrt{-1}$ be a solution.
\item\label{st:len3:out} Output
\[
c_1 := \frac{a+1}{2},\qquad 
c_2 := \frac{a-1}{2}\cdot d_1,\qquad
c_3 := \frac{a-1}{2}\cdot d_2.
\]
\end{enumerate}
\end{alg}

\begin{poc} By assumption, $-1$ is a sum of two squares in~$\K$ and so the norm equation is step~\eqref{st:len3s:norm} is solvable. By a direct computation we obtain that $a$ is the sum of squares of the elements $c_1, c_2, c_3$ constructed in step~\eqref{st:len3:out}. This proves the correctness of the algorithm.
\end{poc}

We can now concentrate on fields of level strictly greater than~$2$ (so either~$4$ or infinity). Observe that the level of a global function field equals the level of its full field of constants and so it cannot exceed~$2$. Thus, if the level is greater than~$2$, the field in question must be a number field.

\begin{alg}\label{alg:len3}
Let~$\K$ be a number field of level $s(K) \neq 1, 2$ and let $a\in \un$ be an element of length $\len(a) = 3$. This algorithm outputs $c_1, c_2, c_3\in \un$ such that $a = c_1^2 + c_2^2 + c_3^2$.
\begin{enumerate}
\item Construct a set $S := \{\gp_1, \dotsc, \gp_2\}$ of finite primes of~$\K$ such that:
  \begin{itemize}
  \item $S$ contains all the dyadic primes of~$\K$;
  \item for every prime~$\gp$, if $\ord_\gp a\odd $, then $\gp\in S$.
  \end{itemize}
\item Construct a basis $\SB = \{ \kappa_1, \dotsc, \kappa_k\}$ of the group $\SingS$ of $S$-singular square classes.
\item Construct a set $S_\infty = \{\gr_1, \dotsc, \gr_r\}$ of all the real places of~$\K$.
\item\label{st:len3:loop} Repeat the following steps: 
  \begin{enumerate}
  \item Build a vector $u = (1, \dotsc, 1)$ of length $r = \card{S_\infty}$.
  \item Build a vector $v = (v_1, \dotsc, v_s)$, indexed by the primes in~$S$, setting:
  \[
  v_i :=
  \begin{cases}
  1 & \text{if }(-1,-1)_\gpi = -1,\\
  0 & \text{if }(-1,-1)_\gpi = 1.
  \end{cases}
  \]
  \item Build a vector $w = (0, \dotsc, 0)$ of length $s = \card{S}$.
  \item\label{st:cij} Construct a matrix $A = (a_{ij})$ with $r$ rows \textup(indexed by the real places of~$\K$\textup) and $k$ columns \textup(indexed by elements in~$\SB$\textup), setting
  \[
  a_{ij} :=
  \begin{cases}
  1 & \text{if }\sgn_{\gr_i} \kappa_j = -1,\\
  0 & \text{if }\sgn_{\gr_i} \kappa_j = 1.
  \end{cases}
  \]
  \item Construct a matrix $B = (b_{ij})$ with $s$ rows and $k$ columns, setting
  \[
  b_{ij} :=
  \begin{cases}
  1 & \text{if }(-1,\kappa_j)_\gpi = -1,\\
  0 & \text{if }(-1,\kappa_j)_\gpi = 1.
  \end{cases}
  \]
  \item Construct a matrix $C = (c_{ij})$ with $s$ rows and $k$ columns, setting
  \[
  c_{ij} :=
  \begin{cases}
  1 & \text{if }(a,\kappa_j)_\gpi = -1,\\
  0 & \text{if }(a,\kappa_j)_\gpi = 1.
  \end{cases}
  \]
  \item Check if the following system of $\FF_2$-linear equations is solvable:
  \begin{equation}\tag{$\clubsuit$}\label{eq:len3}
  \left(\begin{array}{c} A\\\hline B\\\hline C\end{array}\right)
  \cdot \begin{pmatrix} x_1\\\vdots\\ x_k\end{pmatrix}
  = \left(\begin{array}{c} u\\\hline v\\\hline w\end{array}\right).
  \end{equation}
  If it is, denote a solution by $\xi_1, \dotsc, \xi_k$ and exit the loop.
  \item If system~\eqref{eq:len3} is not solvable, find a new prime $\gq\notin S$, replace $S$ by $S\cup \{\gq\}$ and update the basis~$\SB$ accordingly. 
  \item Reiterate the loop.
  \end{enumerate}
\item Set 
\[
b := \kappa_1^{\xi_1}\dotsm \kappa_k^{\xi_k}.
\]
\item\label{st:len3:normL} Set $\L := \K(\sqrt{-1})$ and find a solution $d_1 + d_2\cdot \sqrt{-1}\in \L$ to the norm equation:
\[
\norm{x} = -b.
\]
\item\label{st:len3:normM} Set $M := \K(\sqrt{a})$ and find a solution $d_3 + d_4\cdot \sqrt{a}\in M$ to the norm equation:
\[
\norm[\ext{M}{K}]{y} = b.
\]
\item Output
\[
c_1 := \frac{d_1}{d_4},\qquad 
c_2 := \frac{d_2}{d_4},\qquad 
c_3 := \frac{d_3}{d_4}.
\]
\end{enumerate}
\end{alg}

\begin{poc}
First we need to show that the loop in step~\eqref{st:len3:loop} terminates. That is, that after appending enough primes to the set~$S$, systen~\eqref{eq:len3} eventually becomes solvable. We know that $a$ is a sum of three squares, hence the quadratic form $\form{1,1,1,-a}$ is isotropic over~$\K$. It follows, that there is some $c\in\un$ such that $(-c)$ is represented by the form $\form{1,1}$ and $c$ is represented by $\form{1, -a}$. In other words, the two forms $\form{1, 1, c}$ and $\form{1, -a, -c}$ are both isotropic. In particular, $c$ is totally negative.

Let $T$ be a (finite) set of places of~$\K$, consisting of
\begin{itemize}
\item all archimedean places,
\item all dyadic primes,
\item all non-dyadic primes, at which $a$ has odd valuation.
\end{itemize}
In other words, $T$ is the union of the sets $S$, $S_\infty$ and all the complex places of~$\K$. Then \cite[Lemma~2.1]{LW1992} says that there is a prime $\gq\notin T$ and a $\bigl(T\cup \{\gq \}\bigr)$-singular element $b\in \un$ such that:
{
\renewcommand{\theenumi}{$\textup{C}_{\arabic{enumi}}$}
\begin{enumerate}
\item\label{it:C1} $\sgn_\gr b = \sgn_\gr c$ for every real place~$\gr$ of~$K$,
\item\label{it:C2} $b\equiv c\pmod{\gp^{1+ \ord_\gp 4}}$ for every finite prime $\gp\in T$,
\item $\ord_\gq b = 1$.
\end{enumerate}
}
\noindent Consider the following two quadratic forms:
\[
\q_1 := \form{1, 1, b}\qquad\text{and}\qquad \q_2 := \form{1, -a, -b}
\]
and let $\SB = \{\kappa_1, \dotsc, \kappa_k\}$ be a basis, over~$\FF_2$, of the group $\Sing{T\cup \{\gq\}} = \Sing{S\cup\{\gq\}}$ of $\bigl(S\cup \{\gq\}\bigr)$-singular square classes. Write the square class of~$b$ in a form $\kappa_1^{\xi_1}\dotsm \kappa_k^{\xi_k}$. We will show that the coordinates $\xi_1, \dotsc, \xi_k$ form a solution to~\eqref{eq:len3}, and that every solution makes the corresponding forms~$\q_1$ and~$\q_2$ isotropic. 

Of course, $\q_1$, $\q_2$ are locally isotropic at complex places of~$\K$. Likewise, they are locally isotropic at every finite prime $\gp\notin S\cup \{\gq\}$ by \cite[Corollary~VI.2.5]{Lam2005}, since $\gp$ is non-dyadic and all three coefficients of each form have even valuations at~$\gp$. Therefore, there are only finitely many places of~$\K$ that must be considered. These are precisely the places in $S_\infty\cup S\cup \{\gq\}$. We should deal first with the real places. Observe that $b$, having the same signs as~$c$, is totally negative. Therefore we have
\[
A\cdot \begin{pmatrix} \xi_1\\ \vdots\\ \xi_k \end{pmatrix}
= u.
\]
Conversely, every solution to the above equation gives a totally negative element~$b$. Therefore, the forms $\q_1\otimes \Kr$ and $\q_2\otimes \Kr$ are isotropic for every real place~$\gr$. To prove that the latter one is isotropic, we use the fact that $a$ is totally positive as a sum of squares.

We may now focus on non-archimedean places. Take a prime $\gp\in S$. By \eqref{it:C2} and Local Square Theorem (see e.g. \cite[Theorem~VI.2.19]{Lam2005}) the local squares classes $b\cdot \squares{\Kp}$ and $c\cdot \squares{\Kp}$ coincide. This yields the isometries:
\[
q_1\otimes \Kp\cong \form{1,1,c}\otimes\Kp
\qquad\text{and}\qquad
q_2\otimes \Kp\cong \form{1,-a,-c}\otimes\Kp.
\]
The right hand sides are isotropic and so are the left hand sides. Finally, take the prime~$\gq$ that we appended to~$T$. By the previous part, the two Hilbert symbols:
\[
(-1,-b)_\gp
\qquad\text{and}\qquad
(a,b)_\gp
\]
vanish for every $\gp\neq \gq$. It follows from Hilbert Reciprocity Law (see e.g., \cite[Chapter~VII]{OMeara2000}) that
\[
(-1,-b)_\gq = (a,b)_\gq = 1.
\]
Thus, the forms $\q_1\otimes \Kq$ and $\q_2\otimes \Kq$ are also isotropic. This way, we proved that $\q_1$ and $\q_2$ are locally isotropic everywhere, hence they are isotropic over~$\K$ by the local-global principle (see e.g., \cite[Theorem~VI.3.1]{Lam2005}). 

Toke again a prime $\gp_i\in S\cup \{\gq\}$. We have already proved that $\gq_1\otimes \K_\gpi$ is isotropic. Hence, by \cite[Proposition V.3.22]{Lam2005}, its Hasse invariant $s_\gpi(\q_1)$ equals the Hilbert symbol $(-1-\det \q_1)_\gpi$. This leads to: 
\begin{align*}
(-1)^{v_i}
&= (-1,-1)_\gpi\cdot 1 \\
&= (-1,-1)_\gpi\cdot (1,1)_\gpi\cdot (1,b)_\gpi^2\\
&= (-1,-1)_\gpi\cdot s_\gpi(\q_1)\\
&= (-1,-1)_\gpi\cdot (-1, -\det \q_1)_\gpi\\
&= (-1,-1)_\gpi\cdot (-1, -b)_\gpi\\
&= (-1, b)_\gpi\\
&= \prod_{j\leq k} (-1, \kappa_j)^{\xi_j}_\gpi\\
&= \prod_{j\leq k} (-1)^{\xi_j b_{ij}}.
\end{align*}
Consequently we have
\[
\sum_{j\leq k} \xi_jb_{ij} = v_i.
\]
Conversely, if the above equality holds then $s_\gpi(\q_1) = (-1, -\det \q_1)_\gpi$ and so $\q_1\otimes \K_\gpi$ is isotropic.

Analogously, $\q_2\otimes \K_\gpi$ is isotopic, hence $s_\gpi(\q_2) = (-1, -\det \q_2)_\gpi$. Write
\begin{equation}\label{eq:len3:wi}
\begin{split}
(-1)^{0}
&= (-1,-ab)_\gpi^2\\
&= (-1, -\det\q_2)_\gpi\cdot (-1, -ab)_\gpi \\
&= s_\gpi\q_2\cdot (-1, -ab)_\gpi \\
&= (-a,-b)_\gpi\cdot (-1,-a)_\gpi\cdot (-1,b)_\gpi \\
&= (a,b)_\gpi \\
&= \prod_{j\leq k} (a, \kappa_j)_\gpi^{\xi_j}\\
&= \prod_{j\leq k} (-1)^{\xi_j c_{ij}}.
\end{split}
\end{equation}
This way we obtained the ``$C$-part'' of system~\eqref{eq:len3}:
\[
\sum_{j\leq k} \xi_jc_{ij} = w_i.
\]
Conversely, if the above equality holds for some $\xi_1, \dotsc, \xi_k\in \{0,1\}$ and $b = \kappa_1^{\xi_1}\dotsm \kappa_k^{\xi_k}$, then $s_\gpi \q_2 = (-1, -\det \q_2)_\gpi$ and so $\q_2\otimes \Kpi$ is isotopic.

This way we proved that once we appended~$\gq$ to the initial set~$S$, the system~\eqref{eq:len3} becomes solvable and every solution to this system results in an element $b\in \Sing{S\cup \{\gq\}}$ for which the forms~$\q_1$ and~$\q_2$ are isotropic. 

The existence of the prime~$\gq$ in \cite{LW1992} is proved using Chebotarev density theorem. Although it is not explicit in the statement of \cite[Lemma~2.1]{LW1992}, the set of such primes has positive density. This means that the loop in step~\eqref{st:len3:loop} terminates (see Remark~\ref{rem:prime_q} at the end of the paper).

Now, let $b = \kappa_1^{\xi_1}\dotsm \kappa_k^{\xi_k}$ for some $\xi_1, \dotsc, \xi_k\in \{0,1\}$ forming a solution to~\eqref{eq:len3}. Then, $\q_1 = \form{1, 1, b}$ is isotropic, hence $(-b)$ is a sum of two squares in~$\K$. Therefore, the norm equation in step~\eqref{st:len3:normL} has a solution. Say, 
\begin{equation}\label{eq:len3:d1d2}
-b = \norm{d_1 + d_2\cdot \sqrt{-1}} = d_1^2 + d_2^2.
\end{equation}
Likewise, the norm equation $b = \norm[\ext{M}{K}]{y}$ also has a solution, since $b$ is represented by the form $\form{1, -a}$. Write
\begin{equation}\label{eq:len3:d3d4}
b = \norm[\ext{M}{K}]{d_3 + d_4\cdot \sqrt{a}} = d_3^2 - a\cdot d_4^2.
\end{equation}
Combining \eqref{eq:len3:d1d2} with \eqref{eq:len3:d3d4} we arrive at the sought representation of~$a$ as a sum of three squares:
\[
a = \left( \frac{d_1}{d_4} \right)^2 + \left( \frac{d_2}{d_4} \right)^2 + \left( \frac{d_3}{d_4} \right)^2.\qedhere
\]
\end{poc}

We may now turn our attention to sums of four squares. Recall that the level of a global function field is either~$1$ or~$2$, hence the Pythagoras number of such field is at most~$3$ by \cite[Theorem~XI.5.6]{Lam2005}. Consequently, the $4$-squares problem is meaningful only for number fields. 

\begin{alg}\label{alg:len4}
Let $a\in \un$ be an algebraic number of length $\len(a) = 4$. This algorithm finds $c_1, \dotsc, c_4\in \un$ such that $a = c_1^2 + \dotsb + c_4^2$.
\begin{enumerate}
\item Construct a set $D := \{\gd_1, \dotsc, \gd_d\}$ consisting of the dyadic primes of~$\K$, whose ramification indices and inertia degrees are simultaneously odd, i.e.
\[
e(\ext{\gd_i}{2}) \equiv f(\ext{\gd_i}{2})\odd
\]
for $i\leq d$.
\item Construct a set~$S\supseteq D$ consisting of all the dyadic primes of~$\K$ and all the primes where $a$ has odd valuation.
\item Construct a set $S_\infty = \{\gr_1, \dotsc, \gr_r\}$ of all the real places of~$\K$.
\item\label{st:len4:gihi} For every prime $\gdi\in D$ find elements $g_i, h_i\in \un$ such that
\[
(a, g_i)_\gdi = 1
\qquad\text{and}\qquad
(g_i, h_i)_\gdi = -1.
\]
\item Construct a basis $\SB = \{\kappa_1, \dotsc, \kappa_k\}$ of the group~$\SingS$ of $S$-singular square classes.
\item Repeat the following steps: 
  \begin{enumerate}
  \item Build a vector $u = (1, \dotsc, 1)$ of length $r = \card{S_\infty}$.
  \item Build a vector $v = (1, \dotsc, 1)$ of length $d = \card{D}$.
  \item Build a vector $w = (0, \dotsc, 0)$ of length $s = \card{S}$.
  \item Construct a matrix $A = (a_{ij})$ with $r$ rows \textup(indexed by the real places of~$\K$\textup) and $k$ columns \textup(indexed by elements of~$\SB$\textup), setting:
  \[
  a_{ij} := 
  \begin{cases}
  1 & \text{if }\sgn_{\gr_i}\kappa_j = -1,\\
  0 & \text{if }\sgn_{\gr_i}\kappa_j = 1.
  \end{cases}
  \]
  \item Construct a matrix $B = (b_{ij})$ with $d$ rows \textup(indexed by the primes in~$D$\textup) and $k$ columns, setting:
  \[
  b_{ij} := 
  \begin{cases}
  1 & \text{if }(h_i,\kappa_j)_\gdi = -1\\
  0 & \text{if }(h_i,\kappa_j)_\gdi = 1.
  \end{cases}
  \]
  \item Construct a matrix $C = (c_{ij})$ with $s$ rows \textup(indexed by the primes in~$S$\textup) and $k$ columns, setting
  \[
  c_{ij} := 
  \begin{cases}
  1 & \text{if }(a,\kappa_j)_\gpi = -1\\
  0 & \text{if }(a,\kappa_j)_\gpi = 1.
  \end{cases}
  \]
  \item Check if the following system of $\FF_2$-linear equations is solvable:
  \begin{equation}\tag{$\spadesuit$}\label{eq:len4}
  \left(\begin{array}{c} A\\\hline B\\\hline C \end{array}\right)
  \cdot
  \begin{pmatrix} x_1\\ \vdots\\ x_k \end{pmatrix}
  =
  \left(\begin{array}{c} u\\\hline v\\\hline w \end{array}\right).
  \end{equation}
  If it is, denote a solution by $\xi_1, \dotsc, \xi_k$ and exit the loop.
  \item Otherwise, find a new prime $\gq\notin S$, append it to~$S$ and update the basis~$\SB$. 
  \item Reiterate the loop.
  \end{enumerate}
\item Set $b := \kappa_1^{\xi_1}\dotsm \kappa_m^{\xi_m}$.
\item\label{st:len4:d1d2d3} Execute Algorithm~\ref{alg:len3} to obtain $d_1, d_2, d_3\in \un$ such that
\[
-b = d_1^2 + d_2^2 + d_3^2.
\]
\item\label{st:len4:d4d5} Set $\L := \K\bigl(\sqrt{a}\bigr)$ and find a solution $d_4 + d_5\cdot \sqrt{a}\in L$ to the norm equation
\[
b = \norm{x}.
\]
\item\label{st:len4:output} Output $c_1 := \sfrac{d_1}{d_5}, \dotsc, c_4 := \sfrac{d_4}{d_5}$.
\end{enumerate}
\end{alg}

\begin{poc}
The proof is very similar to that of Algorithm~\ref{alg:len3}. Since $a$ is a sum of four squares, there is some $c\in \un$ such that the two forms $\form{1,1,1,c}$ and $\form{1, -a,-c}$ are isotropic. In particular, $c$ must be totally negative. Again, there is a positive-density set of primes of~$\K$ such that appending to~$S$ any prime~$\gq$ from this set, we may find an $\bigl(S\cup \{\gq\}\bigr)$-singular element~$b$ such that:
{
\renewcommand{\theenumi}{$\textup{C}_{\arabic{enumi}}$}
\begin{enumerate}
\item $\sgn_\gr b = \sgn_\gr c$ for every real place~$\gr$ of~$K$,
\item $b\equiv c\pmod{\gp^{1+ \ord_\gp 4}}$ for every prime $\gp\in S\setminus D$, 
\item $b\equiv g_i\pmod{\gd_i^{1+ \ord_\gdi 4}}$ for the primes $\gd_1, \dotsc, \gd_d\in D$, 
\item $\ord_\gq b = 1$.
\end{enumerate}
}
\noindent By the Local Square Theorem we have
\begin{itemize}
\item $b\cdot \squares{\Kp} = c\cdot \squares{\Kp}$ for every $\gp\in S\setminus D$,
\item $b\cdot \squares{\K_\gdi} = g_i\cdot \squares{\K_\gdi}$ for $\gd_1, \dotsc, \gd_d\in D$.
\end{itemize}
Denote
\[
\q_1 := \form{1,1,1,b}\qquad\text{and}\qquad \q_2 := \form{1,-a,-b}.
\]
As in the proof of correctness of Algorithm~\ref{alg:len3}, we will show that the coordinates $\xi_1, \dotsc, \xi_k$ of~$b$ in $\Sing{S\cup\{\gq\}}$ form a solution to~\eqref{eq:len4}. And every solution to~\eqref{eq:len4} results in an element~$b$ for which $\q_1$ and $\q_2$ are isotropic.

For archimedean places of~$\K$, as well as the non-archimedean ones, but distinct from~$\gq$ and not sitting in~$S$, the arguments are the same as used in the previous proof. Take a prime $\gpi\in S\setminus D$. Then $\q_1\otimes \Kpi$ is isotropic since it contains three units. For the form~$\q_2$ we use the same arguments as in the previous proof to show that $\q_2\otimes \Kpi$ is isotropic, since $b\cdot \squares{\Kpi} = c\cdot \squares{\Kpi}$, and this isotropy is equivalent to the condition
\[
\sum_{j\leq k} \xi_jc_{ij} = w_i.
\]

What really distinguishes sums of four squares from sums of three squares is the local behavior of~$\q_1$ and~$\q_2$ at the dyadic primes $\gd_1, \dotsc, \gd_d\in D$. Fix a prime $\gdi\in D$. We know that $b\cdot \squares{\K_\gdi} = g_i\cdot \squares{\K_\gdi}$ and $(g_i, h_i)_\gdi = -1$. Therefore we have
\[
(-1)^{v_i}
= -1
= (h_i, b)_\gdi
= \prod_{j\leq k} (h_i, \kappa_j)^{\xi_j}_\gdi
= \prod_{j\leq k} (-1)^{\xi_jb_{ij}}.
\]
This means that the following $\FF_2$-linear condition holds
\[
\sum_{j\leq k} \xi_j b_{ij} = v_i.
\]
On the other hand, if $\xi_1, \dotsc, \xi_k\in \{0,1\}$ satisfy the following equality, and $b = \kappa_1^{\xi_1}\dotsm \kappa_k^{\xi_k}$, then $(h_i, b) = -1$ and so $b$ is not a local square at~$\gdi$. Consequently the form $\q_1\otimes \Kdi$ is isotropic by \cite[Corollary VI.2.15]{Lam2005}.

Similar arguments apply to the form $\q_2\otimes \Kdi$. By assumption, $e(\ext{\gdi}{2})$ and $f(\ext{\gdi}{2})$ are odd. Hence \cite[Example~XI.2.4]{Lam2005} asserts that $s(\Kdi) = 4$ and consequently $(-1,-1)_\gdi = -1$. Using the facts that $b\cdot \Kdi = g_i\cdot \Kdi$ and $(a, g_i)_\gdi = 1$ we may write
\begin{align*}
s_\gdi\q_2
&= (-a, -b)_\gdi \\
&= (-1, -ab)_\gdi\cdot (a, b)_\gdi \\
&= (-1, -\det \q_2)_\gdi\cdot (a, g_i)_\gdi \\
&= (-1, -\det \q_2).
\end{align*}
Therefore, $\q_2\otimes \Kdi$ is isotropic by \cite[Proposition~V.3.22]{Lam2005}. Repeating the same calculations that we did in Eq.~\eqref{eq:len3:wi} we see that
\[
\sum_{j\leq k} \xi_jc_{ij} = w_i.
\]
Conversely, if the above equality holds for some $\xi_1, \dotsc, \xi_k\in \{0,1\}$ and $b = \kappa_1^{\xi_1}\dotsm \kappa_k^{\xi_k}$, then $s_\gdi\q_2 = (-1, -\det \q_2)_\gdi$ and so $\q_2\otimes \Kdi$ is isotropic, as desired.

Finally, take the prime~$\gq$. The form $\q_1\otimes \Kq$ is isotropic because it contains three units and $\gq$ is non-dyadic. On the other hand, using Hilbert Reciprocity Law in the same way as we did in the proof of correctness of Algorithm~\ref{alg:len3}, we show that $\q_2\otimes \Kq$ is isotropic, as well.

All in all, we proved that the coordinates of the square class of~$b$ in the group~$\Sing{S\cup\{\gq\}}$, treated as a linear space over~$\FF_2$, form a solution to~\eqref{eq:len4}. And that for every such solution the corresponding forms~$\q_1$ and~$\q_2$ are locally isotropic at every completion of~$\K$, hence the local-global principle says that they are isotropic over~$\K$. This implies that the algorithm terminates. The rest of the proof is fully analogous to the previous one, where we dealt with sums of three squares.
\end{poc}

\begin{rem}\label{rem:prime_q}
In Algorithms~\ref{alg:len3} and~\ref{alg:len4} we keep enlarging the set~$S$ until we happen to come across a prime~$\gq$ from a certain positive-density set. There are two basic strategies to do it, with different worst-case scenarios and average behaviors. One method is to select primes of~$\K$ at random. The probability that we will find a ``good'' prime~$\gq$ is given by the density of this set, and so is positive. However, in the worst case scenario, we may hypothetically keep picking the primes that are constantly not in this set. Hence, while this method is better-behaved in practice, it does not warrant that the algorithms terminate in finite time. In theory it may take indefinitely long to find a correct prime~$\gq$. The other approach uses a purely deterministic exhaustive search. Start with a (rational) prime $p=3$ and try all the primes of~$\K$ extending~$p$. Then go to the next prime number, take the primes of~$\K$ extending it, and so on. Although slower in practice, this procedure guaranties that the two algorithms terminate in finite time.
\end{rem}

\begin{rem}
In step~\eqref{st:len4:gihi} of Algorithm~\ref{alg:len4}, for every $\gdi\in D$ with $i\leq \card{D}$, we need to find elements $g_i, h_i\in\K$ such that 
\[
(a, g_i)_\gdi = 1
\qquad\text{and}\qquad 
(g_i, h_i)_\gdi = -1. 
\]
The element~$g_i$ is not used again in the algorithm. It is needed only to find~$h_i$ and can be discarded afterwards. (On the other hand it is very important in the proof of correctness of the algorithm.) 

Here again, a fast and efficient method is a probabilistic one. Let $n_i := (\Kdi : \QQ_2)$. If $a \notin \squares{\Kdi}$ then the probability that a random element $g_i\in K$ satisfies the first of the above two conditions is~$\sfrac12$. However, the probability that $g_i$ itself is a local square (and so the second condition can never be satisfied) is $2^{-n_i-2}$. Thus, the probability that we pick a ``good'' $g_i$ at random equals $(\sfrac12 - \sfrac{1}{2^{n_i+2}}) \geq \sfrac14$. On the other hand, if $a\in \squares{\Kdi}$, then \emph{every} $g_i\in K$ satisfies the first condition. Nevertheless, $g_i$ still cannot be a local square at~$\gdi$, for the second condition to be satisfiable. Therefore, the probability of finding a ``good'' $g_i$ at random is $(1 - \sfrac{1}{2^{n_i+2}}) \geq \sfrac34$. In both cases, once we selected~$g_i$, an element~$h_i$ satisfying the second condition can be found by a random search with probability~$\sfrac12$.
\end{rem}

\section{Conclusion}
In this paper we present an algorithmic method for finding a representation of an element of a global field as a sum of squares with the minimal number of summands. The described algorithms have been implemented by the author using computer algebra system Magma \cite{BCP1997} as a part of package CQF \cite{Koprowski2020}.

\bibliographystyle{acm} 
\bibliography{sos}
\end{document}